\documentclass[a4paper,12pt]{amsart}

\usepackage{amsmath,amsthm,amssymb}
\usepackage{bbm}
\usepackage{mathrsfs}
\usepackage{enumerate}
\usepackage{graphicx}
\usepackage{caption} 
\usepackage{subcaption} 
\usepackage{thmtools}
\usepackage{thm-restate}

\newcommand{\R}{\mathbb{R}}

\newcommand{\del}{\partial }

\newtheorem{thm}{Theorem}

\theoremstyle{definition}

\newtheorem{question}[thm]{Question}
\theoremstyle{remark}
\newtheorem{rem}[thm]{Remark}
\begin{document}

\title{Trisections of 3--manifold bundles over $S^1$}
\author{Dale Koenig}
\address{Topology and Geometry of Manifolds Unit, Okinawa Institute of Science and Technology, Okinawa, Japan}
\email{dale.koenig@oist.jp}
\urladdr{https://groups.oist.jp/manifolds/dale-koenig}
\date{}

\begin{abstract}
Let $X$ be a bundle over $S^1$ with fiber a 3--manifold $M$ and with monodromy $\varphi$.  Gay and Kirby showed that if $\varphi$ fixes a genus $g$ Heegaard splitting of $M$ then $X$ has a genus $6g+1$ trisection.  Genus $3g+1$ trisections have been found in certain special cases, such as the case where $\varphi$ is trivial, and it is known that trisections of genus lower than $3g+1$ cannot exist in general.  We generalize these results to prove that there exists a trisection of genus $3g+1$ whenever $\varphi$ fixes a genus $g$ Heegaard surface of $M$.  This means that $\varphi$ can be nontrivial, and can preserve or switch the two handlebodies of the Heegaard splitting.  We additionally describe an algorithm to draw a diagram for such a trisection given a Heegaard diagram for $M$ and a description of $\varphi$.
\end{abstract}

\maketitle

\section{Introduction}
A trisection of a closed, smooth, four dimensional manifold $X$ is a decomposition $X = X_1 \cup X_2 \cup X_3$ such that each $X_i$ is a 4--dimensional handlebody $\natural^{k_i} S^1 \times B^3$ and each pairwise intersection $X_i \cap X_j$, $i \not = j$ is a 3--dimensional handlebody $\natural^g S^1 \times D^2$.  The triple intersection $X_1 \cap X_2 \cap X_3$ is a closed surface of genus $g$, and induces a Heegaard splitting of each $\del X_i$.  Since each $X_i$ is a four dimensional handlebody, $\del X_i$ is diffeomorphic to $\#^{k_i} S^1 \times S^2$ for some nonnegative integer $k_i$.  If $k_1 = k_2 = k_3$ then the trisection is called \emph{balanced} and we say that it is a $(g;k_1)$--trisection of $X$.  Otherwise it is called \emph{unbalanced} and we say that it is a $(g;k_1,k_2,k_3)$--trisection of $X$.  Here $g$ denotes the \emph{genus} of the trisection.  The triple intersection surface, together with marked cut systems indicating the boundaries of disks in each of the three $X_i$ is called a \emph{diagram} of the trisection.  A diagram for a trisection completely determines the trisection itself.  For more details we refer the reader to the original paper by Gay and Kirby \cite{GayKirby1}. 

\begin{figure}[h]
\centering
\includegraphics[height=1.6in]{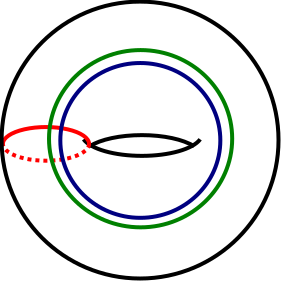}
\caption{One of the three possible genus 1 trisection diagrams of $S^4$ up to diffeomorphism of the diagram surface.  The other two are obtained by recoloring.}
\label{s3g1}
\end{figure}

We also will need to use the notion of (unbalanced) stabilization and destabilization of trisections.  These are described in detail in \cite{MeierSchirmerZupan1}.  We only need the fact that an unbalanced stabilization affects the trisection diagram by connect summing with one of the three genus one trisections of $S^4$, one of which is shown in Figure \ref{s3g1}.  We can similarly view destabilization in terms of diagrams.  Suppose we have a trisection diagram given by surface $\Sigma$ and cut systems $\alpha,\beta,\gamma$.  Suppose there is some simple closed curve $C$ in $\Sigma$ that cuts off a genus 1 trisection diagram of $S^4$.  Then compressing the surface along $C$ and removing the genus 1 diagram for $S^4$ results in a destabilization of the original trisection. 

Let $M$ be a smooth 3--manifold with a Heegaard splitting $H_1 \cup_\Sigma H_2$.  Let $\varphi \colon M \to M$ be a self--diffeomorphism of $M$, possibly isotopic to the identity.  We say that $\varphi$ \emph{preserves} the Heegaard splitting if $\varphi(H_1) = H_1$ and $\varphi(H_2)=H_2$.  We say that $\varphi$ \emph{flips} the Heegaard splitting if $\varphi(H_1) = H_2$ and $\varphi(H_2)=H_1$.  In both cases $\varphi(\Sigma) = \Sigma$; that is, $\varphi$ sends the Heegaard surface to itself.  For any choice of $\varphi$, both $H_1 \cup_\Sigma H_2$ and $\varphi(H_1) \cup_{\varphi(\Sigma)} \varphi(H_2)$ are Heegaard splittings of $M$, and if $\Sigma'$ is a stabilization of $\Sigma$ then $\varphi(\Sigma')$ is a stabilization of $\varphi(\Sigma)$.  It follows from the Reidemeister--Singer theorem that after enough stabilizations we get a Heegaard splitting that is preserved by $\varphi$, up to isotopy \cite{Reidemeister1}\cite{Singer1}.  For the remainder of the paper, we fix notation to assume that $X$ is always a 4-manifold constructed as above with $\varphi$ either preserving or flipping the Heegaard splitting $H_1 \cup_\Sigma H_2$ of the fiber $M$.  We henceforth use $g$ to denote the genus of $\Sigma$.  For notational convenience we let $H$ denote an abstract genus $g$ handlebody, so $H$ is in the same diffeomorphism class as both $H_1$ and $H_2$.

In the case where $\varphi$ preserves a genus $g$ Heegaard splitting of $M$, Gay and Kirby showed that $X$ has a $(6g+1;2g+1)$--trisection \cite{GayKirby1}.  In this paper we use similar techniques to improve on this result.  We prove the following:

\begin{thm}
\label{mainresult}
Let $X$ be a bundle over $S^1$ with fiber a closed 3--manifold $M$.  If the monodromy $\varphi$ preserves or flips a genus $g$ Heegaard splitting of $M$, then $X$ has a $(3g+1;g+1)$--trisection.
\end{thm}

We will prove the case where $\varphi$ preserves the splitting and the case where it flips the splitting separately.  Perhaps unintuitively, the case where $\varphi$ flips the splitting is easier, so we will prove the theorem in this case first.

\begin{rem}
Since we care about how $\varphi$ acts on the Heegaard surface and not just how it acts on $M$, we can have interesting monodromy even when $\varphi$ is isotopic to the identity.  For example, if $\varphi$ is the identity but $H_1$ is isotopic to $H_2$, then $\varphi$ will act on $\Sigma$ nontrivially, sending each cut system for $H_1$ to a cut system for $H_2$.
\end{rem}

After constructing these trisections, we will describe how to construct diagrams for them.  In the case where $\varphi$ flips the Heegaard splitting, we derive the diagram directly.  In the case where $\varphi$ preserves the Heegaard splitting, we will first construct a diagram for a $(4g+1;2g+1,g+1,g+1)$ unbalanced trisection of $X$ and show that we can destabilize down to diagram of the desired $(3g+1;g+1)$ balanced trisection.

There is a correspondence between trisections of a 4--manifold $X$ and certain handle decompositions of $X$.  A trisection can be turned into a handle decomposition of $X$ where $X_1$ corresponds to the 0--handle and 1--handles, $X_2$ to the 2--handles, and $X_3$ to the 3--handles and 4--handle.  Since the rank of the fundamental group cannot be lower than the number of 1--handles, $k_1$ gives a lower bound for the rank of $\pi_1(X)$.  By the symmetry of $X_1,X_2,X_3$, $k_2$ and $k_3$ do as well.  Suppose $M = \#^g S^1 \times S^2$ and $\varphi$ is trivial.  Then $\varphi$ fixes the standard genus $g$ Heegaard splitting of $M$, so there exists a $(3g+1;g+1)$--trisection of $X$.  The rank of $\pi_1(X)$ is equal to $g+1$ and each $X_i$ is genus $g+1$,  so each $X_i$ is the lowest genus possible.  Since the genera of the $X_i$ together with the Euler characteristic of $M$ determine the genus of the trisection, there exist no trisections of lower genus.  It follows that Theorem \ref{mainresult} cannot be improved in general, although it is possible that there exist lower genus trisections for specific choices of $M$ and $\varphi$.

Much of the work in this paper is covered in \cite{koenig1}, which shows that there exists a $(3g+1;g+1)$ trisecction when $\varphi$ flips a genus $g$ Heegaard splitting of $M$, or in the case where $\varphi$ is trivial.  In this paper we expand to the case where $\varphi$ is nontrivial but preserves the Heegaard splitting, and provide a more detailed explanation of how to construct diagrams from these trisections.

\subsection*{Acknowledgements} Special thanks go to Jeffrey Meier for his comments on the previous draft of this paper, and to the author's advisor Abigail Thompson for her encouragement and many useful conversations.

\section{The Main Theorem}
We split into the two cases.

\smallskip
\noindent\textit{Case 1: $\varphi$ flips the Heegaard splitting}\ \ 

\begin{figure}[h]
\centering
\includegraphics[height=1.6in]{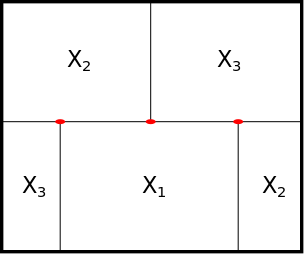}
\caption{$I \times M$ is depicted as a rectangle with $I$ as the horizontal axis and $M$ as the vertical axis.  The middle horizontal line represents $I \times \Sigma$, so the copy of $M$ represented by each vertical slice is cut by the middle line into copies of the handlebodies $H_1$ and $H_2$.}
\label{flippablebreakdown}
\end{figure}

We view $X$ as being constructed by gluing together the two ends of $I \times M$.  Cut $I \times M$ into regions as depicted in Figure \ref{flippablebreakdown}.  Note that there are two regions labelled $X_2$ and two labelled $X_3$.  The top regions are copies of $I \times H_2$ and the bottom regions are copies of $I \times H_1$.  Since $\varphi$ flips the Heegaard splitting, the two $X_2$ regions are combined into a single $I \times H$ region, and the two $X_3$ regions are combined into another $I \times H$ region.  There is only a single $X_1$ region that is already diffeomorphic to $I \times H$.  Therefore, this already divides $X$ into three 4 dimensional handlebodies.  However, it is not yet a trisection because the pairwise intersections $X_i \cap X_j$ are not 3 dimensional handlebodies.  In fact,  they are not even connected.

To fix this we make a modification at each horizontal segment in the two dimensional picture lying in some $X_i \cap X_j$.  Note that there are three such segments, since the two segments shown in $X_2 \cap X_3$ were glued together by $\varphi$.   Consider the horizontal segment between $X_1$ and $X_2$, representing a copy of $I \times \Sigma$ in $X_1 \cap X_2$.  Let $B$ be a small ball in $M$ intersecting $\Sigma$ in a single disk.  Remove $I \times B$ from $X_1 \cup X_2$ and assign it to $X_3$ instead.  Thus, we have run a tube from $X_3$ along the horizontal segment, pushing apart $X_1$ and $X_2$ to connect the two sides of the $X_3$ region.  Repeat the process two more times, running a tube from $X_1$ along the horizontal segment lying in $X_2 \cap X_3$ and a tube from $X_2$ along the horizontal segment lying in $X_1 \cap X_3$.  The choices of $B$ for different tubes must chosen to be disjoint, for example by choosing small neighborhoods of different points on $\Sigma$.  Note that in any generic vertical slice, the boundary of the 3--ball $B$ consists of a disk bordering the $X_i$ region below $B$ and a disk bordering the $X_j$ region above $B$.

We now show that $X_1 \cap X_2$ is now a genus $g+1$ handlebody.  The argument for the other two pairwise intersections is identical.  $X_1 \cap X_2$ consists of the following parts:

\begin{itemize}
\item A copy of $I \times D^2$ coming from the $X_1$ tube through $X_2 \cap X_3$.
\item A copy of $I \times (\Sigma - D^2)$ coming from the horizontal segment between $X_1$ and $X_2$.  This is the product of an interval and a punctured genus $g$ surface, so it is a genus $2g$ handlebody.  
\item A copy of $I \times D^2$ coming from the $X_2$ tube through $X_1 \cap X_3$.
\item A copy of $H$ coming from the vertical segment between the $X_1$ and $X_2$ labelled regions.  This is a genus $g$ handlebody.
\end{itemize}

Therefore, $X_1 \cap X_2$ consists of a genus $2g$ and a genus $g$ handlebody connected by two different three dimensional 1--handles.  The result is a genus $3g+1$ handlebody.  Since each $X_i$ is a copy of $I \times H$ with the tube acting as an additional 1--handle, each is a genus $g+1$ four dimensional handlebody.  Since each pairwise intersection is a handlebody and $X_1 \cap X_2 \cap X_3$ is the boundary of each pairwise intersection, the triple intersection is necessarily a closed surface.  It follows that we have indeed constructed a $(3g+1;g+1)$--trisection of $X$.

\smallskip
\noindent\textit{Case 2: $\varphi$ preserves the Heegaard splitting}\ \ 

\begin{figure}[h]
\centering
\includegraphics[height=1.6in]{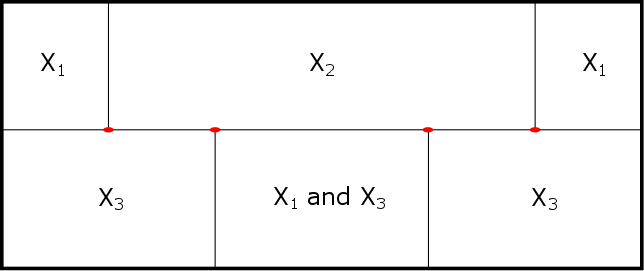}
\caption{$I \times M$ is depicted as a rectangle with $I$ as the horizontal axis and $M$ as the vertical axis.  The middle horizontal line represents $I \times \Sigma$, so the copy of $M$ represented by each vertical slice is cut by the middle line into copies of the handlebodies $H_1$ and $H_2$.}
\label{preservebreakdown}
\end{figure}

We again begin by cutting up $I \times M$, this time using Figure \ref{preservebreakdown} as a guide.  Note that since $\varphi$ preserves the Heegaard splitting, the two $X_1$ regions and the two $X_3$ regions are glued together to create a single $X_1$ region and a single $X_2$ region.  The region labelled ``$X_1$ and $X_3$" needs some explanation.  This region is a copy of $I \times H_1$.  Decompose $H_1$ into two three dimensional balls, $U$ and $V$, that intersect in $g+1$ disks.  In this region, $I \times U$ is assigned to $X_1$ and $I \times V$ is assigned to $X_3$.  We then attach three tubes, the first connecting the two sides of the $X_2$ region through the horizontal line belonging to $X_1 \cap X_3$.  Attach the other two tubes to $X_1$, connecting the region labelled $X_1$ with $I \times U$ through the two horizontal lines belonging in $X_2 \cap X_3$.  As before, we denote the tubes by $I \times B$ where the choices of $B$ for each tube are disjoint.  Each $X_i$ is the union of $I \times H$ with one extra four dimensional 1--handle attaching $\{0\} \times H$ and $\{1\} \times H$, so each $X_i$ is a genus $g+1$ handlebody.  It remains to check that each pairwise intersection is a genus $3g+1$ three dimensional handlebody.  We check this for each pair.

$X_1 \cap X_2$ breaks down as follows:

\begin{itemize}
\item Two copies of $I \times D^2$ coming from the tubes connecting the $X_1$ labelled region with the ``$X_1$ and $X_3$" region.  There is one tube for each of the two horizontal lines in $X_2 \cap X_3$.
\item One copy of $I \times D^2$ coming from the $X_2$ tube through the horizontal line in $X_1 \cap X_3$.
\item Two copies of $H_2$ coming from the vertical lines in $X_1 \cap X_2$.
\item The product of an interval and a $g+1$ punctured sphere, coming from the intersection of $I \times U$ and $X_2$.  This piece is diffeomorphic to a genus $g$ handlebody.
\end{itemize}

Each of the two copies of $H_2$ contributes $g$ to the genus.  The product of the interval and a $g+1$ punctured sphere contributes another $g$ to the genus.  This gives three disconnected pieces, which are then connected up by the three $I \times D^2$ 1--handles to give a single connected piece.  Since only two 1--handles would be needed to connect the pieces up, the third adds an additional 1 to the genus.  It follows that the result is indeed a genus $3g+1$ three dimensional handlebody.

$X_1 \cap X_3$ breaks down as:

\begin{itemize}
\item A copy of $I \times (\Sigma - D^2)$ coming from the horizontal segment between $X_1$ and $X_3$.  This is the product of an interval and a punctured genus $g$ surface, so it is a genus $2g$ handlebody.  
\item A copy of $I \times D^2$ from each of the two $X_1$ tubes through the corresponding horizontal $X_2 \cap X_3$ segments.
\item A copy of $U$ (a three dimensional ball) lying at each of the two vertical boundaries between the $X_3$ and ``$X_1$ and $X_3$" regions.
\item $g+1$ copies of $I \times D^2$ coming from the $g+1$ disk components of intersection between $U$ and $V$ in the ``$X_1$ and $X_3$" region.
\end{itemize}

The union of the two copies of $U$ with the $g+1$ copies of $I \times D^2$ is a genus $g$ three dimensional handlebody.  Indeed, it is the result of connecting two balls by $g+1$ 1--handles.  This piece and the genus $2g$ handlebody coming from the horizontal $X_1 \cap X_3$ segment are connected together with two different 1--handles, adding 1 more to the genus.   The result is then a genus $3g+1$ handlebody as expected.

Lastly, $X_2 \cap X_3$ breaks down as:
\begin{itemize}
\item Two copies of $I \times (\Sigma - D^2)$ coming from the two horizontal lines between $X_2$ and $X_3$.
\item A copy of $I \times D^2$ from the tube of $X_2$ through the $X_1 \cap X_3$ horizontal segment.
\item The product of an interval and a $g+1$ punctured sphere, coming from the intersection of $I \times V$ and $X_2$.  This piece is diffeomorphic to a genus $g$ handlebody.
\end{itemize}

Note that $\Sigma - D^2$ can be obtained from the $g+1$ punctured sphere by attaching $g$ additional 1--handles.  Therefore, the union of the $I \times \{\mathrm{punctured~sphere}\}$ piece and the two $I \times (\Sigma - D^2)$ pieces is diffeomorphic to a genus $3g$ handlebody.  The copy of $I \times D^2$ then adds one more to the genus, for $3g+1$ total.

Since each pairwise intersection is indeed a genus $3g+1$ handlebody, we have successfully constructed a $(3g+1;g+1)$--trisection of $X$.  We constructed genus $3g+1$ trisections in both the cases where $\varphi$ preserves and where it flips the Heegaard splitting, so the main theorem is proven.

\section{The diagrams}
Again we split into the two separate cases.  We proceed from the previous descriptions of the trisections, finding cut systems for each of the three pairwise intersections $X_i \cap X_j$.  Again, it is easier to start with the case where $\varphi$ flips the Heegaard splitting.

\smallskip
\noindent\textit{Case 1: $\varphi$ flips the Heegaard splitting}\ \ 

\begin{figure}[h]
\centering
\includegraphics[height=3in]{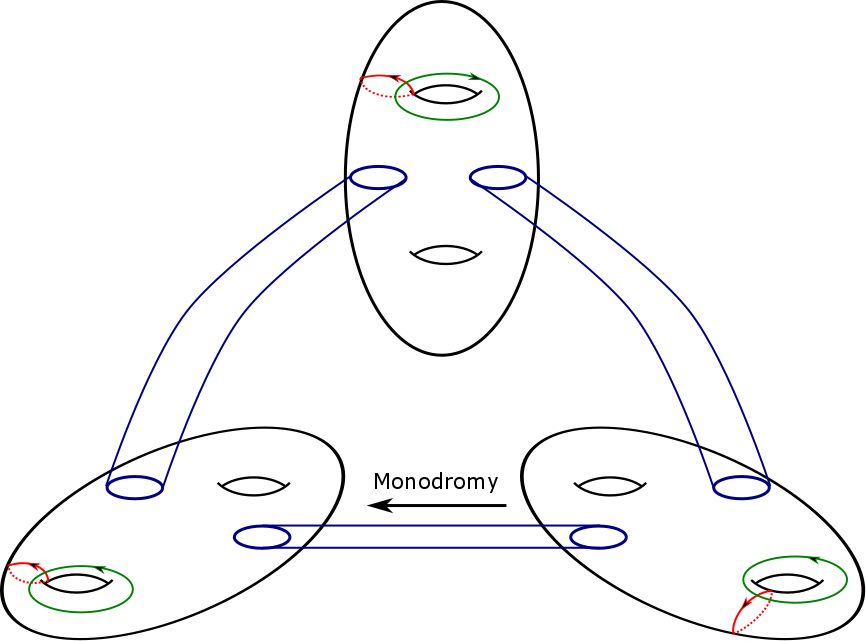}
\caption{  The top copy of $\Sigma$ has been reflected about $\theta = \pi/4$.  The red and green oriented curves in each copy of $\Sigma$ therefore are the same pair of curves slid along $I \times \Sigma$.  Passing through the gluing monodromy, however, may send the green and red curve to something else.  We view the monodromy as all happening between the bottom two copies of $\Sigma$. }
\label{embeddingsigmas}
\end{figure}

We begin by describing the triple intersection surface.  The reader may have noticed in Figure \ref{flippablebreakdown} that there was a red dot lying at each of the three points lying in the intersection of the three regions.  At the beginning of the construction, before attaching the tubes to $X_1,X_2$ and $X_3$ along the horizontal $I \times \Sigma$ segments, $X_1 \cap X_2 \cap X_3$ consisted of three copies of $\Sigma$, one lying at each red dot.  After adding the tubes, $X_1 \cap X_2 \cap X_3$ consists of three twice punctured copies of $\Sigma$ together with a a tube connecting each pair.  This results in a genus $3g+1$ surface.  We construct the diagram keeping this in mind, drawing three copies of $\Sigma$ and connecting them with tubes.  Additionally, we can choose freely how to embed these copies of $\Sigma$, so we embed them in a way that makes the curves that bound disks look especially simple.   Specifically, we draw the three copies of $\Sigma$ along three axes lying in the plane and evenly distributed in $\theta$, and then reflect $\Sigma_2$ locally in $\theta$.  That is, if $\theta_2$ is the value of $\theta$ at which we drew $\Sigma_2$, reflect $\Sigma_2$ about the ray $\theta = \theta_2$.  See Figure \ref{embeddingsigmas}.  

\begin{figure}[h]
\centering
\includegraphics[height=3in]{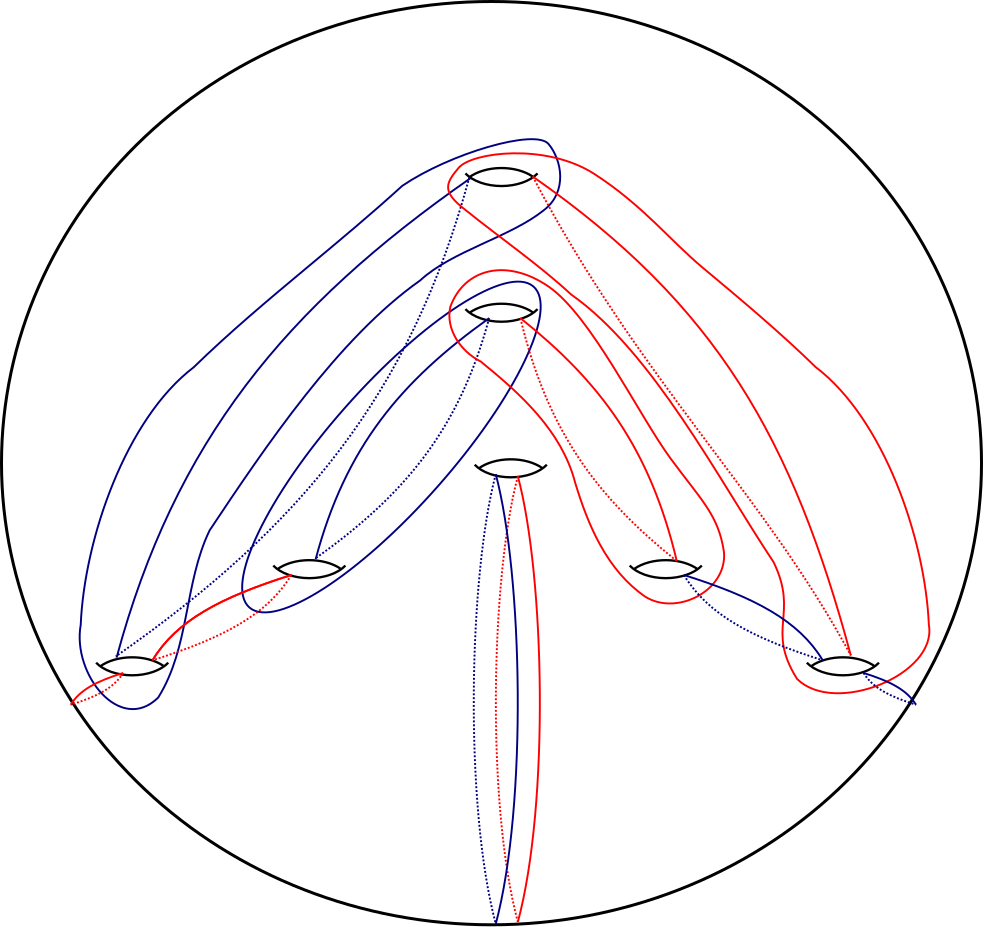}
\caption{The first two sets of curves when $g=2$.}
\label{firsttwocolors}
\end{figure}

When drawing diagrams we draw the first copy of $\Sigma$ on the bottom left, the second on top, and the third in the bottom right.  Denote these by $\Sigma_1, \Sigma_2$ and $\Sigma_3$ respectively.  Let $\theta_{ij}$ denote the angle midway between where we draw $\Sigma_i$ and $\Sigma_j$.  With the embeddings chosen for $\Sigma_i$, nontrivial curves that look symmetric about $\theta_{12}$ and avoid $\Sigma_3$ bound disks in one of the $I \times \{\mathrm{punctured~surface}\}$ segments.  The same is true for curves that look symmetric about $\theta_{23}$ and avoid $\Sigma_1$.  The curves bounding disks in the last copy of $I \times \{\mathrm{punctured~surface}\}$ avoid $\Sigma_2$, but will look more complicated in general because they must pass through the monodromy.   The algorithm to draw a diagram for the trisection is as follows:

\begin{enumerate}
\item Draw the genus $3g+1$ surface as copies of $\Sigma_1$, $\Sigma_2$, and $\Sigma_3$ with tubes connecting each pair.
\item Draw $2g$ disjoint nonisotopic blue curves symmetric about $\theta_{12}$ and $2g$ disjoint nonisotopic red curves symmetric about $\theta_{23}$.  The blue curves neccessarily avoid $\Sigma_3$ and similarly the red curves necessarily avoid $\Sigma_1$.  These are the curves bounding product disks in two of the $I \times \{\mathrm{punctured~surface}\}$ segments.
\item Draw $g$ red curves as meridians of $\Sigma_1$, and $g$ blue curves as meridians of $\Sigma_3$.
\item Draw a parallel red and blue curve as meridians of the tube connecting  $\Sigma_1$ and $\Sigma_3$.  At this point, the diagram should look like Figure \ref{firsttwocolors}.
\item On $\Sigma_2$ draw in green the set of curves that bound disks in $H_2$.  That is, draw a Heegaard diagram for $M$ in green.  Keep in mind the special embedding of $\Sigma_2$, so the diagram will reflected in $\theta$ relative to a diagram drawn on $\Sigma_1$ or $\Sigma_3$.
\item Draw a green curve as a meridian of the tube connecting $\Sigma_2$ with either $\Sigma_1$ or $\Sigma_3$.  The two choices are easily seen to be handle slide equivalent.
\item Let $c_1,\cdots,c_{2g}$ be a set of disjoint nonparallel arcs in $\Sigma_3$ with their boundaries lying on the puncture of $\Sigma_3$ that is tubed to $\Sigma_1$.  Then $\varphi(c_1),\cdots,\varphi(c_{2g})$ is a set of disjoint nonparallel arcs in $\Sigma_1$ with their boundaries lying on the puncture of $\Sigma_1$ that is tubed to $\Sigma_3$.  Draw $2g$ green curves by connecting up each $c_i$ with $\varphi(c_i)$ using two arcs in the tube between $\Sigma_1$ and $\Sigma_3$.  If $a,b$ are the two endpoints of $c_i$, we connect $c_i$ and $\varphi(c_i)$ such that $a$ is connected to $\varphi(a)$ and $b$ to $\varphi(b)$.  In other words, we connect each arc to its image with the opposite orientation.  These are the curves bounding product disks in the final  $I \times \{\mathrm{punctured~surface}\}$ segment.
\end{enumerate}

\begin{figure}[h]
\centering
\includegraphics[height=3in]{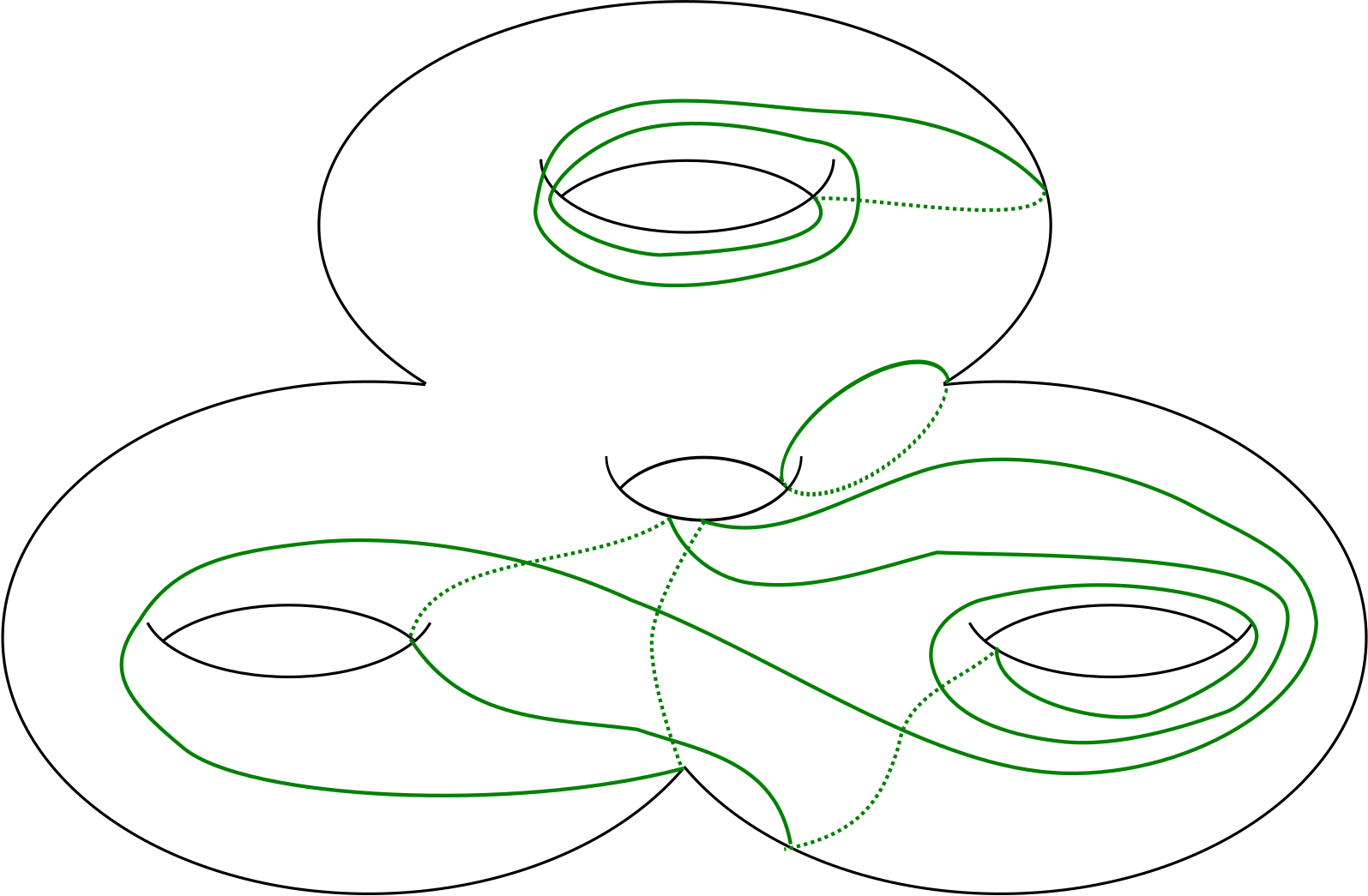}
\caption{The third set of curves when $M =\R P^3$ and $\varphi$ flips $H_1$ and $H_2$.  Restricting the bottom two green curves to $\Sigma_3$ gives a complicated pair of arcs that turn into trivial looking arcs when $\varphi$ is applied to the surface.   Connecting the former arcs with the latter, keeping orientations of the arcs in mind, produces the two green curves shown.}
\label{greenex}
\end{figure}

Note that red and blue curves found to bound product disks in the second step are found in the same way as the green curves in the last step, except that there is no monodromy between $\Sigma_2$ and the other two surfaces, so the image of the arcs is determined purely by the embeddings of the $\Sigma_i$.  The choices of embeddings imply that the image of an arc in $\Sigma_1$ will be its image in $\Sigma_2$ reflected across $\theta_{12}$, so the closed curves obtained look symmetric about $\theta_{12}$.  The same is true when obtaining the curves symmetric about $\theta_{23}$.  For an example of how the green curves might look, Figure \ref{greenex} shows the set of green curves for the case where $M$ is $\R P^3$ and the monodromy isotopic to the identity but flipping the Heegaard splitting, sending $H_1$ to $H_2$ and vice versa.

\smallskip
\noindent\textit{Case 2: $\varphi$ preserves the Heegaard splitting}\ \ 

\begin{figure}[h]
\centering
\includegraphics[height=1.6in]{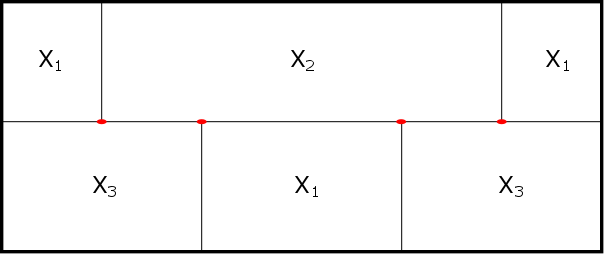}
\caption{$I \times M$ is depicted as a rectangle with $I$ as the horizontal axis and $M$ as the vertical axis.  The middle horizontal line represents $I \times \Sigma$, so the copy of $M$ represented by each vertical slice is cut by the middle line into copies of the handlebodies $H_1$ and $H_2$.}
\label{unbalancedbreakdown}
\end{figure}

We begin by constructing the genus $4g+1$ unbalanced trisection shown in Figure \ref{unbalancedbreakdown} by attaching tubes through the areas indicated by horizontal lines as in the previous constructions.  We will construct a diagram for this trisection and then later destabilize the diagram to get a genus $3g+1$ diagram.  We draw triple intersection surface for this diagram as four copies of $\Sigma$ connected by tubes.  We call the four copies $\Sigma_1, \Sigma_2, \Sigma_3$, and $\Sigma_4$ where $\Sigma_i$ of cyclically adjacent indices are conneccted by tubes.  We draw them in the plane at evenly distributed $\theta$ values as before, by drawing one with all the holes at the same $\theta$, and rotating it through the plane to get the other three copies.  In the examples we draw $\Sigma_1$ at the top right, and then draw $\Sigma_2, \Sigma_3, \Sigma_4$ clockwise.  We reflect each of $\Sigma_2$ and $\Sigma_4$ locally in $\theta$ as we did with $\Sigma_2$ in the case where $\varphi$ flips the Heegaard splitting.  This ensures that we can choose simple looking curves bounding product disks in the copies of $I \times (\Sigma - D^2)$ except between $\Sigma_1$ and $\Sigma_4$ where we assume the monodromy map occurs.  Again let $\theta_{ij}$ be the value of $\theta$ midway between where we draw $\Sigma_i$ and $\Sigma_j$.  The algorithm to draw the diagram on this surface is then:

\begin{enumerate}
\item Draw the genus $4g+1$ surface as copies of $\Sigma_i$, $i = 1,2,3,4$ with adjacent indices connected by tubes.
\item Draw $2g$ disjoint nonisotopic blue curves symmetric about $\theta_{12}$ and avoiding $\Sigma_3$ and $\Sigma_4$, and another $2g$ blue curves symmetric about $\theta_{34}$ and avoiding $\Sigma_1$ and $\Sigma_2$.
\item Draw a single blue curve as a meridian of the tube connecting $\Sigma_1$ and $\Sigma_4$.
\item Draw $2g$ disjoint nonisotopic green curves symmetric about $\theta_{23}$ and avoiding $\Sigma_1$ and $\Sigma_4$.
\item Draw a single red and green curve each around either the tube connecting $\Sigma_1$ and $\Sigma_2$ or the tube connecting $\Sigma_3$ and $\Sigma_4$.  The two choices are slide equivalent.
\item Draw a Heegaard diagram for $M$ in green on $\Sigma_1$.
\item Draw another Heegaard diagram for $M$ in green on $\Sigma_4$, but remember the special embedding of $\Sigma_4$.  Thus, this diagram will appear to be reflected in $\theta$ compared to the diagram drawn in the last step.
\item Draw a meridian system for each of $\Sigma_2$ and $\Sigma_3$ in red, for $2g$ red curves total.
\item Let $c_1,\cdots,c_{2g}$ be a set of disjoint nonparallel arcs in $\Sigma_1$ with their boundaries lying on the puncture of $\Sigma_1$ that is tubed to $\Sigma_4$.  Then connect each $c_i$ to its image $\varphi(c_i)$ in $\Sigma_4$ with the opposite orientation.  Draw these $2g$ curves in red.  Remember that the image will need to be drawn reflected in $\theta$ because of the embedding chosen for $\Sigma_4$.  If $\varphi$ acts trivially on $\Sigma$ these curves will look symmetric about $\theta_{41}$.
\end{enumerate}

\begin{figure}[h]
\centering
\includegraphics[height=3in]{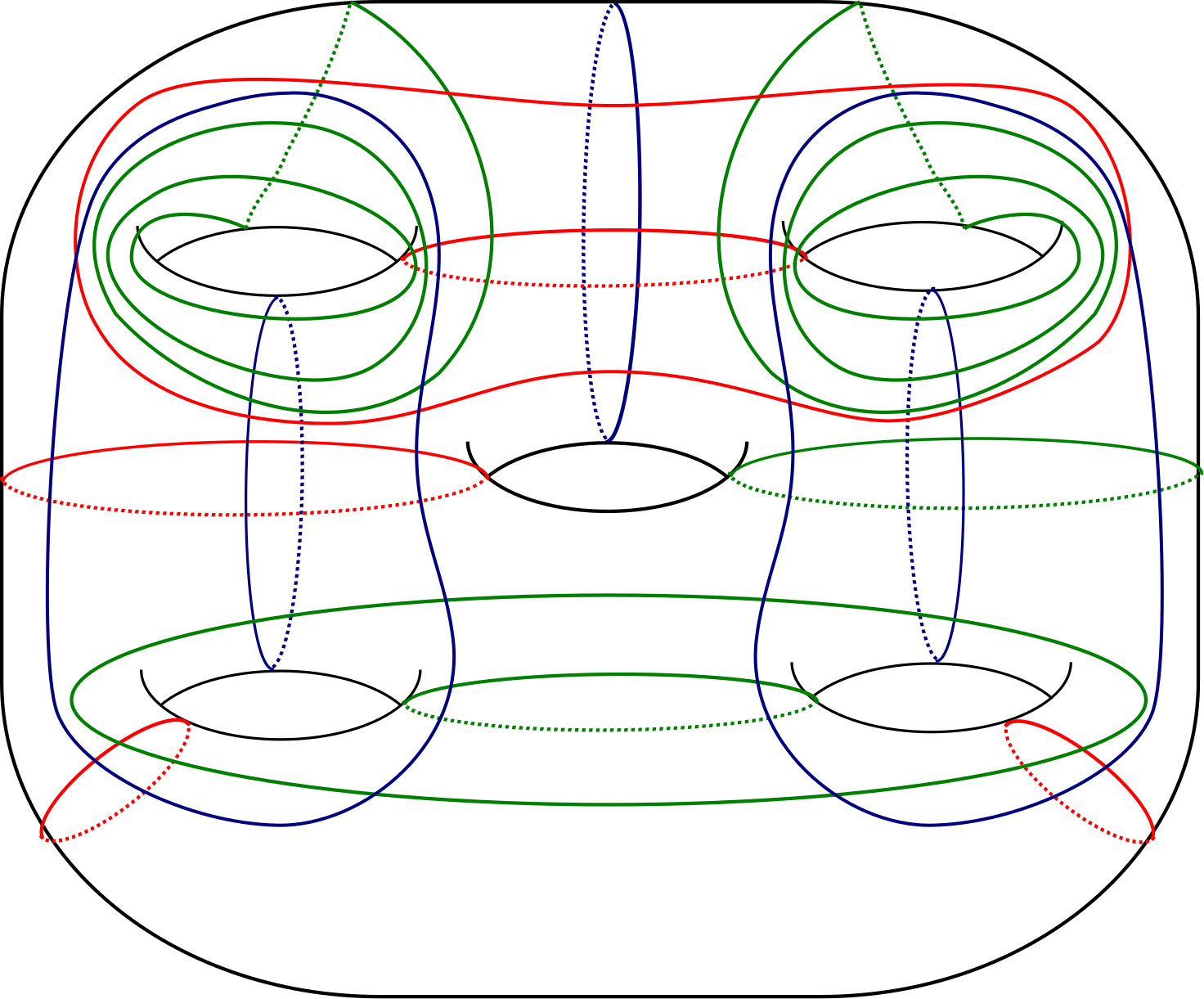}
\caption{A diagram for the genus $5$ trisection where $M = L(3,1)$ and $\varphi$ is the identity. }
\label{unbalancedex}
\end{figure}

An example is shown in Figure \ref{unbalancedex}.

\begin{figure}[h]
\centering
\includegraphics[height=3in]{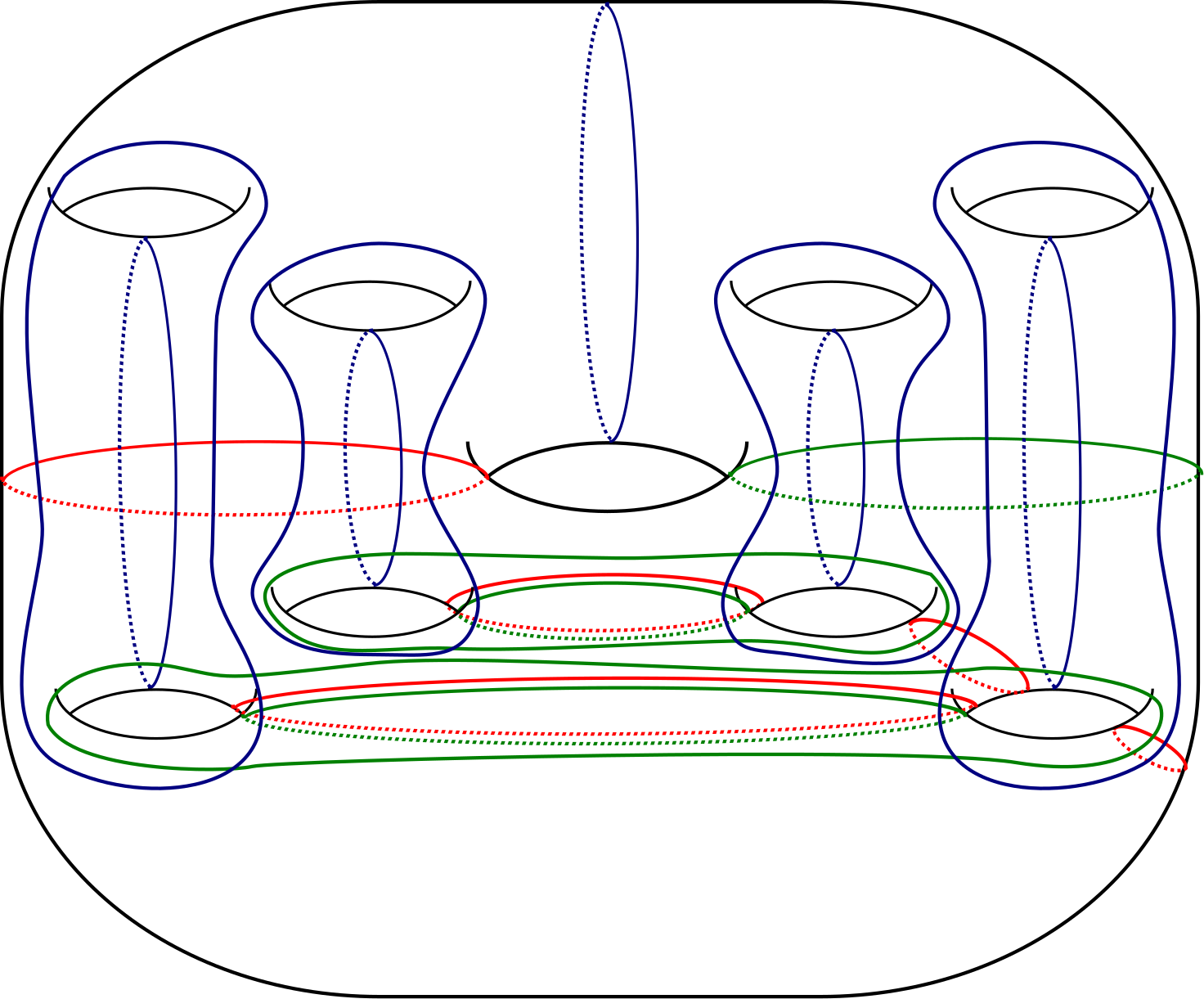}
\caption{A partial trisection diagram for for finding destabilizations in the case where $g=2$.  Note that additional curves that would be required to complete the diagram will not interfere with the destabilizations.}
\label{destabpartial}
\end{figure}

\begin{figure}[h]
\centering
\includegraphics[height=3in]{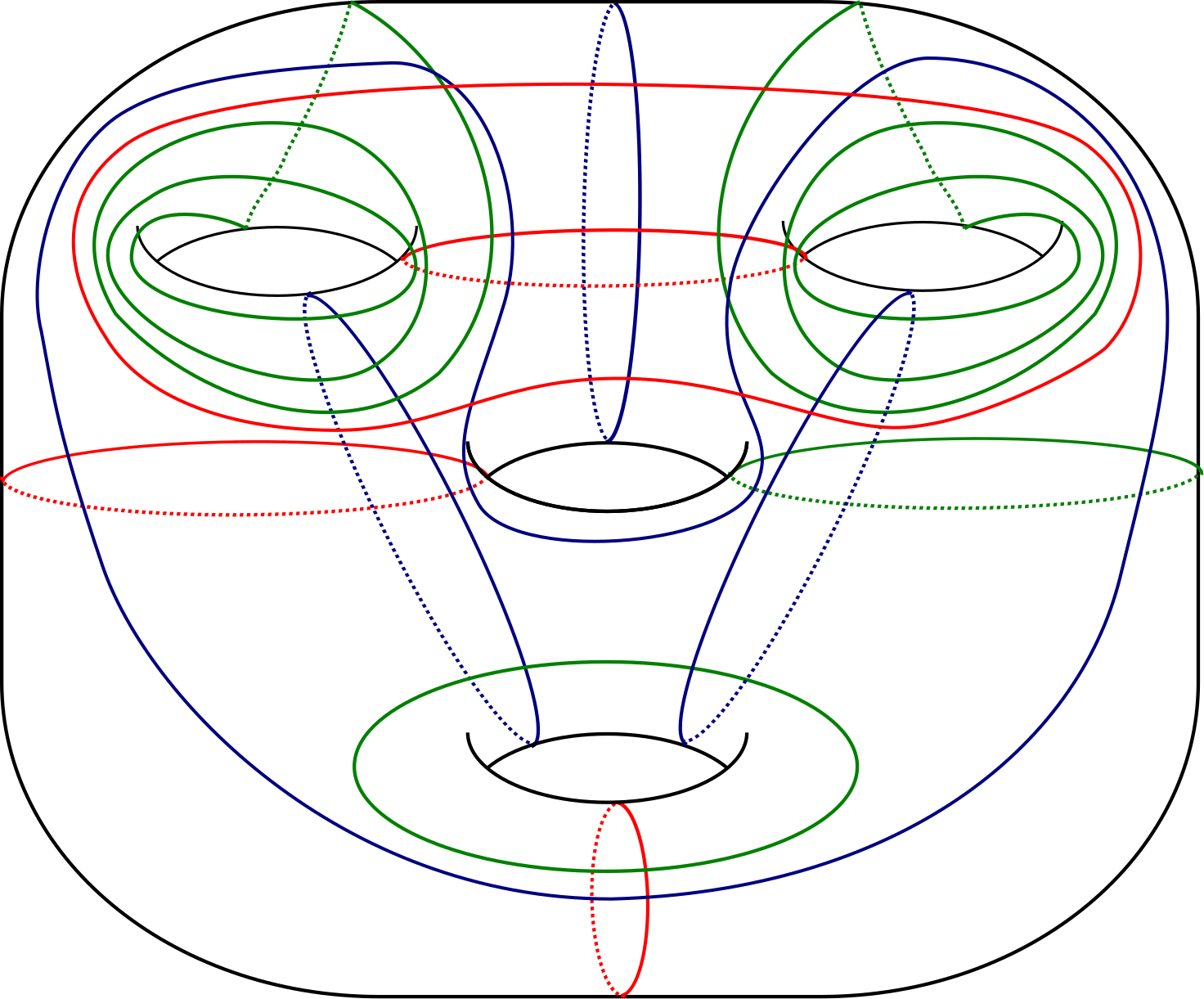}
\caption{A diagram for the genus 4 trisection where $M = L(3,1)$ and $\varphi$ is the identity. }
\label{ls12destab}
\end{figure}

We now show that this diagram can always be destabilized $g$ times.  First, choose the symmetric blue and green curves to look as in Figure \ref{destabpartial}, consisting of curves that either look like meridians in the chosen embedding or that lie on the near side of the surface and loop around one of those meridians.   This is already true in Figure \ref{unbalancedex}.  Next we can slide the red meridian curves lying on $\Sigma_2$ to be parallel to the $g$ green curves that look like meridians as we have done in Figure \ref{destabpartial}.  At this point we have $g$ pairs of parallel red and green curves, and the isotopy class of this curve intersects two of the blue curves in one point each.  See Figure \ref{destabpartial} for a partial Heegaard diagram.  For each pair, we can slide one of these blue curves over the other, after which we have a parallel red and green curve that each intersect the set of blue curves in a single point.  This implies that our diagram is a ($g$ times) unbalanced stabilization, and the destabilization can be achieved by compressing along the isotopy classes of each parallel red and green curve, and then deleting the now trivial red and green curves together with the blue curve that they intersect.  See Figure \ref{ls12destab} for the diagram resulting from destabilizing the trisection in Figure \ref{unbalancedex}.  Therefore, we can always get a genus $3g+1$ diagram by constructing a genus $4g+1$ diagram using the algorithm above, and then performing these $g$ destabilizations.  An analysis of how these destabilizations affect the trisecction itself shows that these destabilizations actually produce the genus $3g+1$ trisection constructed in the main theorem.

\section{Further Questions}
We end by leaving a few unanswered questions for the interested reader.
\begin{question}
Is there some choice of $M$ and $\varphi$ such that there exists a trisection of genus lower than $3g+1$, where $g$ is the minimal genus of a Heegard splitting of $M$?
\end{question}

We next ask if we get the same in cases where we could apply both methods to get two trisections of $X$.

\begin{question}
Suppose that $\varphi$ can be chosen up to composition with an isotopy to either preserve or flip a Heegaard splitting of $M$.   Are the trisections produced in these two cases equivalent?
\label{unequal1}
\end{question}

A more specific and possibly easier version of the question above is the following.
\begin{question}
Suppose $X = S^1 \times M$ and $\Sigma$ is a Heegaard surface for a flippable Heegaard splitting of $M$.  That is, $H_1$ and $H_2$ are isotopic.  Then $\varphi$ could be chosen to be the identity, preserving the Heegaard splitting, or it could be chosen to be nontrivial but isotopic to the identity and to flip the Heegaard splitting.  These two choices give two trisections of $X$.  Are these two trisections always equivalent?
\label{unequal2}
\end{question}

Lastly, we note that the construction in the case where $\varphi$ preserves the Heegaard splitting creates the trisection in an assymetric manner, where the $X_i$ play different roles.  This leads to the following questions:
\begin{question}
Is there a more symmetric construction of a $(3g+1;g+1)$--trisection in the case where $\varphi$ preserves the Heegaard splitting?  Is there a symmetric or nearly symmetric diagram?
\end{question}

\begin{question}
Suppose $\varphi$ preserves the Heegaard splitting and we use the methods of this paper to construct a $(3g+1;g+1)$--trisection.  Does rotating the indicies of the $X_i$ always produce an equivalent trisection?
\label{unequal3}
\end{question}

An answer of ``no" to one of questions \ref{unequal1}, \ref{unequal2}, or \ref{unequal3} would be especially interesting as there are, as of the time of this paper, no known examples of 4--manifolds $X$ with two distinct balanced trisections of the same genus.

\bibliographystyle{amsplain}
\bibliography{4dbundles}{}

%

%

\end{document}